\documentclass[11pt]{amsart}
\usepackage{amsmath, amssymb, amsthm,color}
\usepackage{amssymb,longtable}
\usepackage{amsfonts}
\usepackage{setspace}
\usepackage{amsmath} 
\usepackage{hyperref} 
\usepackage{latexsym}
\usepackage{array}
\usepackage{amssymb}
\theoremstyle{plain}

\numberwithin{equation}{section}
\newtheorem{thm}{Theorem}[section]

\newtheorem{lemm}[thm]{Lemma}

\theoremstyle{remark}

\newcommand{\bq}{\begin{equation}}
\newcommand{\eq}{\end{equation}}

\begin{document}

\title{Note on The Euler equations in $C^k$ spaces}
\author{Tarek M. Elgindi and Nader Masmoudi}
\date{}
\maketitle

\begin{abstract}
In this note, using the ideas from our recent article \cite{EM}, we prove strong ill-posedness for the 2D Euler equations in $C^k$ spaces. This note provides a significantly shorter proof of many of the main results in \cite{BLi2}. In the case $k>1$ we show the existence of initial data for which the $kth$ derivative of the velocity field develops a logarithmic singularity immediately. The strong ill-posedness covers $C^{k-1,1}$ spaces as well. The ill-posedness comes from the pressure term in the Euler equation. We formulate the equation for $D^k u$ as:

$$\partial_t D^k u=D^{k+1} p + l.o.t.$$

and then use the non-locality of the map $u\rightarrow p$ to get the ill-posedness. The real difficulty comes in how to deal with the "l.o.t." terms which can be handled by special commutator estimates.

\end{abstract}

\section{introduction}
In this note we give a short proof of the strong ill-posedness of the Euler equations for incompressible flow in $C^k$ spaces. Our proof works in the whole space case as well as the case of periodic boundary conditions and the bounded domain case. We use the non-locality of the map $u\rightarrow p$ to get growth in $C^k$ spaces (since Riesz transforms are unbounded on these spaces!). The proof uses classical commutator estimates which estimate the commutation of a Riesz transform (or a Calder\'on Zygmund operator) and composition by a bi-Lipschitz map. 

We prove the following theorem. 

\begin{thm}
The Euler equations are strongly ill-posed in $C^k$ spaces for $k\geq 1$. In other words, for every $\epsilon>0$ there exists initial data $u_0\in C^k$ such that the unique solution, $u(t)$, of the Euler equations with initial data $u_0$ leaves $C^k$ immediately. 

\end{thm}

We note that very recently Bourgain and Li have proven the same result as above \cite{BLi2}. In their paper the authors cite our work \cite{EM} and note that it did not include ill-posedness $C^k$ spaces for $k>1$. In this note we show that the work in \cite{EM} can easily be extended to the case $k>1.$ Moreover we make clear the fact that strong ill-posedness in $C^k$ can be proven quite easily only using commutator estimates without having to rely upon very intricate constructions. 

\section{The proof}
\begin{proof}
Note that it suffices to consider the two dimensional Euler equations. 
Now consider the equation for $ D^k u$ which means $k$ spatial derivatives of $u,$ which we take to be a vector of many components. 

$\nabla D^{k-1} u$ satisfies the following equation:

$$\partial_t \nabla D^{k-1} u + u \cdot \nabla D^k u +\sum_{j,l}^k Q(D^j u, D^l u) +D^{k-1} D^2 p=0 $$

Recall that $$\Delta p= det (\nabla u)$$ so that $$(D^2 p)_{ij} =\big (R_{i}R_{j} det(\nabla u) \big )_{ij}$$

We rely upon the following commutator estimate in $L^p:$

\begin{lemm}
Let $\Phi$ be a bi-Lipschitz measure preserving map. Let $K=\max\{|\Phi-Id|_{Lip}, |\Phi^{-1}-Id|_{Lip} \}.$ Let $R$ be a composition of Riesz transforms. Define the following commutator: $$[R,\Phi] \omega= R(\omega\circ\Phi)-R(\omega)\circ \Phi,$$ for $\omega\in L^p.$

Then, $$\|[R,\Phi]\|_{L^p\rightarrow L^p}\leq c_{p} K.$$ Moreover, $c_p\leq c p$ as $p\rightarrow \infty$.

\end{lemm}

Now we recall the Lagrangian flow 

$$\dot{\Phi}(x,t)= u(\Phi(t,x),t)$$
$$\Phi(x,0)=x.$$

Because $u$ is divergence free, $\Phi$ is measure preserving. Furthermore, $\Phi(x,-t) = \Phi^{-1} (x,t).$ 

Now, we may write $$\Phi(x,t) = x +\int_{0}^{t} u(\Phi(x,\tau),\tau)d\tau.$$

Thus, $$\Phi(\cdot,t) - I = \int_{0}^{t} u(\Phi(\cdot,\tau),\tau)d\tau.$$

Consequently, $$|\Phi-I|_{\text{Lip}} \leq t | u|_{\text{Lip}}|\Phi|_{Lip}$$ and similarly for $\Phi^{-1} (\cdot,t)= \Phi(\cdot,-t).$ 

Furthermore, by Gronwall's lemma, $$|\Phi|_{\text{Lip}} \leq \exp (t | u|_{\text{Lip}}).$$

In particular, 

\begin{equation} |\Phi-I|_{\text{Lip}} \leq t | u|_{\text{Lip}}\exp(t|u|_{\text{Lip}}). \end{equation}

In particular, if $u$ is Lipschitz then the Lagrangian flow-map is controlled.  

Now, assume that $u$ remains Lipschitz (which in the case $k>0$ is trivial). 

Then we have:

$$ \partial_t(D^k u \circ \Phi) +\sum_{j,l}^k Q(D^j u, D^l u)\circ \Phi + \Big ( R_{i}R_{j}D^{k-1}det (\nabla u)\Big ) \circ \Phi=0 $$

Now assume that the solution remains in $C^k.$

Now call $R_i R_j:=R$ for short. 

$$ \partial_t(D^k u \circ \Phi) +\sum_{j,l}^k Q(D^j u, D^l u)\circ \Phi + \Big ( R D^{k-1}det (\nabla u)\Big ) \circ \Phi=0 $$

$$ \partial_t(D^k u \circ \Phi) +\sum_{j,l}^k Q(D^j u, D^l u)\circ \Phi + R \Big (  D^{k-1}det (\nabla u) \circ \Phi\Big )= [R,\Phi]D^{k-1}det (\nabla u) $$

We are going to take a special initial data $u_0\in C^k.$ Now, suppose that the solution remains in $C^k.$ Note that we can always solve the 2D Euler equations in $W^{k,p}$ and that the solution will be unique when $k>1$ (in the case k=1 we will also get a unique solution the class of velocity fields with bounded curl). Assume that the solution stays in $C^k$ with $|u|_{C^k}\leq M.$

Then we can solve the Euler equations formally using the Duhamel formula:

$$|D^k u\circ \Phi|_{L^p} \geq |D^k u_0+tR(D^{k-1} det(\nabla u_0)|_{L^p}-\Big | \int_{0}^t e^{R(t-s)}[R,\Phi]D^{k-1} det(\nabla u) -Q(D^j u, D^l u)\circ \Phi ds  \Big |_{L^p}$$

Now suppose that we construct compactly supported initial data $u_0$ such that $|RD^{k-1} det(\nabla u_0)|_{L^p}\geq cp$ as $p\rightarrow \infty.$

Then we see that 

$$|D^k u\circ \Phi|_{L^p} \geq tcp -C- t(1+tp) |[R,\Phi]|_{L^p\rightarrow L^p} C(M)+tC(M).$$

This is because $\|R\|_{L^p\rightarrow L^p}\lesssim p$ as $p\rightarrow \infty.$
Now using the commutator estimate in Lemma 2.1 we get:

$$|D^k u\circ \Phi|_{L^p} \geq tcp -C- t(1+tp) tp \tilde{C}(M)+tC(M)$$

Now take $t$ very small depending upon $M$ then we get:

$$|D^k u\circ \Phi|_{L^p} \geq t\tilde{c}p$$

for $t$ small and all $p$ large. 

This contradicts the fact that $|D^k u|_{L^\infty}$ remains bounded.
 
\end{proof}

Now, of course we relied upon the existence of an initial data $u_0$ such that $$|RD^{k-1} det(\nabla u_0)|_{L^p}\geq cp.$$

Constructing such initial data is not too difficult. We will copy the construction from our work \cite{EM} below in the case $k=1$. The higher order cases are similar.

We are interested in showing that for some $i,j$ and for some divergence free $u,$ with $\nabla u \in L^\infty,$ $D^2 p= R_{i}R_{j} \text{det}(\nabla u)$ has a logarithmic singularity. Once that is shown, Lemma 8.2 will follow by a regularization argument.  

Take a harmonic polynomial, $Q,$ which is homogeneous of degree 4. In the two-dimensional case, we can take $$Q(x,y):=x^4+y^4-6x^2y^2,$$
$$\Delta P =0.$$

Define $$G(x,y):= Q(x,y) Log (x^2 +y^2).$$

Notice that \begin{equation} \partial_{i}\partial_{j}\Delta G\in L^\infty (B_{1}(0)), i,j\in\{1,2\}.\end{equation}

Notice, on the other hand, that \begin{equation} \partial_{xxyy} G=-24Log(x^2+y^2)+H(x,y),\end{equation} with $H\in L^\infty(B_{1}(0)).$ In particular, $\partial_{xxyy}G$ has a logarithmic singularity at the origin--and the same can be said about $\partial_{xxxx}G$ and $\partial_{yyyy}G.$ 

\vspace{3mm}

Define $\tilde{u}=\nabla^\perp \Delta G.$ Then, by 8.12,  $\nabla \tilde{u} \in L^\infty(B_{1}(0)).$ Moreover, by definition, $$R_{i}R_{j} \nabla \tilde{u}=\nabla\nabla^\perp\partial_{ij}G.$$

Thus, for example, $R_{1}R_2 \nabla\tilde{u}_{1x}=\partial_{xxyy}G$ has a logarithmic singularity in $B_{1}(0).$ Unfortunately, we are interested in showing that $R_{i}R_{j} \text{det}(\nabla u)$ has a logarithmic singularity for some $i,j,$ not $R_{i}R_{j}\nabla u.$ To rectify this, we choose $$u=\delta \nabla^\perp\Delta (\chi G)+\eta(2y\chi+y^2 \partial_y\chi,y\partial_x\chi),$$ where $\eta, \delta$ are small parameters which will be determined and $\chi$ is a smooth cut-off function with: 

$$\chi=1 \, \, \text{on} \, \, B_1(0),$$
$$\chi=0 \,\, \text{on} \, \, B_2(0)^c.$$ 

Note that $u$ is divergence free and $$u\equiv \delta\nabla^\perp \Delta G +  \eta(y,0) \,\, \text{on} \,\, B_{1}(0). $$

Therefore, $$\nabla u= \delta \left[ {\begin{array}{cc} -\partial_{xy} \Delta G & -\partial_{yy}\Delta G \\ \partial_{xx}\Delta G & \partial_{xy}\Delta G \\\end{array} } \right] + \eta   \left[ {\begin{array}{cc}0 & 1 \\0 & 0 \\\end{array} } \right] .$$

In particular, $$\text{det}(\nabla u) = \eta\delta\partial_{xx}\Delta G + \delta^2 J(x,y),$$ where $J$ is a bounded on $B_{1}(0)$.

Now consider $R_2 R_2 \text{det}(\nabla u):$

$$R_2 R_2 \text{det}(\nabla u)= \eta \delta \partial_{xxyy}G +\delta^2 R_{2} R_2 J.$$
Now, by (8.13), we have $$ R_2 R_2 \text{det}(\nabla u)= \eta \delta (-24Log(x^2+y^2)+H(x,y)) +\delta^2 R_{2} R_2 J, $$

with $H$ and $J$ bounded. Now, recall that $R_2R_2$ maps $L^\infty$ to BMO and that any BMO function can have at most a logarithmic singularity. 

Thus,  

$$|R_2 R_2 \text{det}(\nabla u)|\geq 24 \eta \delta Log(x^2+y^2) -C\delta^2 Log(x^2+y^2) -|H(x,y)|.$$

Choose $\delta < < \eta$ and we see that, near $(0,0)$  $$|R_2 R_2 \text{det}(\nabla u)|\geq 12 \delta^2 Log(x^2+y^2).$$

Taking $\delta\leq C$ small enough, we see that $|\nabla u| \leq 1$ but $|R_2 R_2 \text{det} (\nabla u)| \geq c Log(x^2 +y^2), $ for some small $c.$ 

This completes the construction.

\end{document}